\documentclass[11pt]{article}

\usepackage{cite}

\usepackage{amsmath}
\usepackage{amsfonts}
\usepackage{amssymb}
\usepackage{amsthm}
\usepackage{enumitem}
\usepackage{mathtools}
\usepackage{setspace}
\usepackage{etoolbox}

\newtheoremstyle{mystyle}
{11pt}                          
{11pt}                          
{}                                      
{}                                      
{\bfseries}                     
{}                                      
{5.5pt}                         
{}                                      

\theoremstyle{mystyle}
\newtheorem{theorem}{Theorem}[section]

\newtheorem{corollary}[theorem]{Corollary}

\renewenvironment{proof}[1][Proof.]{\vspace{-16.5pt} \begin{trivlist}
        \item[\hskip \labelsep {\bfseries #1}]}{\qed \end{trivlist}}

\setitemize{noitemsep, topsep=5.5pt, parsep=5.5pt, partopsep=0pt}

\allowdisplaybreaks
\appto\normalsize{
        \abovedisplayskip=5.5pt plus 2pt minus 2pt
        \belowdisplayskip=5.5pt plus 2pt minus 2pt
        \abovedisplayshortskip=5.5pt plus 2pt minus 2pt
        \belowdisplayshortskip=5.5pt plus 2pt minus 2pt}
\appto\small{
        \abovedisplayskip=5.5pt plus 2pt minus 2pt
        \belowdisplayskip=5.5pt plus 2pt minus 2pt
        \abovedisplayshortskip=5.5pt plus 2pt minus 2pt
        \belowdisplayshortskip=5.5pt plus 2pt minus 2pt}

\newcommand{\gap}{\vspace{11pt}}

\newcommand{\diag}{\operatorname{diag}}
\newcommand{\tr}{\operatorname{tr}}

\newcommand{\calc}{\mathcal{C}}
\newcommand{\calk}{\mathcal{K}}
\newcommand{\cala}{\mathcal{A}}
\newcommand{\call}{\mathcal{L}}
\newcommand{\E}{\mathcal{E}}

\newcommand{\Sn}{\mathcal{S}^n}

\newcommand{\Snp}{\mathcal{S}_+^n}

\newcommand{\I}{\mbox{I}}
\newcommand{\calz}{\mathcal{Z}}
\newcommand{\calh}{\mathcal{H}}

\newcommand{\lng}{\langle}
\newcommand{\rng}{\rangle}
\newcommand{\lf}{\left}
\newcommand{\rg}{\right}

\newcommand{\re}{\mathbb R}
\newcommand{\rn}{\mathbb R^n}
\newcommand{\an}{\mathcal A^n}
\newcommand{\sn}{\mathcal S^n}
\newcommand{\snp}{\mathcal S^n_+}

\newcommand{\nn}{\mathbb R^{n\times n}}

\newcommand{\tp}{^\top}

\DeclareMathOperator{\copos}{copos}
\DeclareMathOperator{\compos}{compos}

\title{\bf    
The cone of $\calz$-transformations on the second order cone}
\author{S\'andor Z. N\'emeth\\School of Mathematics, University of Birmingham, Watson Building, Edgbaston\\Birmingham B15 2TT, United
Kingdom\\email: s.nemeth@bham.ac.uk
\and\\
        M. Seetharama Gowda\\Department of Mathematics and Statistics, University of Maryland, Baltimore County\\1000 Hilltop
	Circle, Baltimore, Maryland 21250, U.S.A\\email: gowda@umbc.edu
}

\date{}

\begin{document}

\maketitle

\begin{abstract}
	In this paper, we describe the structural properties of the cone of $\calz$-transformations on the second order cone 
 in terms of  
the semidefinite cone and copositive/completely positive cones induced by the second order cone and its boundary. In particular, we 
describe its dual as a slice of the semidefinite cone as well as a slice of the completely positive cone 
of the second order cone. This  provides an example of an instance where a conic linear program on a 
completely positive cone is reduced to a  problem on the semidefinite cone. 
\end{abstract}

\vspace{1cm}
\\

\newpage

\section{Introduction}

Given a proper cone $\calk$ in a finite dimensional real Hilbert space $(\calh,\lng\cdot,\cdot\rng)$, 
a linear transformation 
$A:\calh\rightarrow \calh$ is said to be a {\it $\calz$-transformation on $\calk$} if 
$$\Big [\,x\in \calk, y\in \calk^*,\,\mbox{and}\,\,\langle x,y\rangle=0\,\Big ]\Rightarrow \langle Ax,y\rangle \leq 0,$$
where $\calk^*$ denotes the dual of $\calk$ in $\calh$. 
Such transformations appear in various areas including economics, dynamical systems, optimization, see e.g., 
\cite{BermanPlemons1994,BermanNeumannStern1989,GowdaTao2009,FanTaoRavindran2017} and the references therein. 
When $\calh$ is $\rn$ and $\calk$ is the nonnegative orthant, $\calz$-transformations become $\calz$-matrices, which are square matrices with nonpositive off-diagonal entries. 


\smallskip

The set $\calz(\calk)$ of all $\calz$-transformations on $\calk$ is a closed convex cone in the space of all (bounded) linear transformations on $\calh$. Given their appearance and importance in various areas, describing/characterizing elements of 
$\calz(\calk)$ and its interior, boundary, dual, etc., is of interest. 
An early result of Schneider and Vidyasagar \cite{SchneiderVidyasagar1970} asserts that $A$ is a $\calz$-transformation on $\calk$ if and only if
$e^{-tA}(\calk)\subseteq \calk$ for all $t\geq 0$; consequently,
\begin{equation}\label{sv result}
\calz(\calk)=\overline{\re\, \I-\pi(\calk)},
\end{equation}
where $\pi(K)$ denotes the set of all linear transformations that leave $\calk$ invariant, $\I$ denotes the identity transformation, and overline denotes the closure. To see another description of $\calz(\calk)$, 
let   
 LL$(\calk):=\calz(\calk)\cap -\calz(\calk)$ denote the lineality space of $\calz(\calk)$, the elements of which are called Lyapunov-like transformations. 
Then the inclusions
$$\re \,\I-\pi(\calk)\,\subseteq \,\mbox{LL}(\calk)-\pi(\calk)\,\subseteq\, \calz(\calk)=\overline{\re \,\I-\pi(\calk)}$$
imply that $$\calz(\calk)=\overline{ \mbox{LL}(\calk)-\pi(\calk)}.$$
As the cones $\calz(\calk)$, $\pi(\calk)$, and LL$(\calk)$  are generally difficult to describe for an arbitrary
proper cone  $\calk$, we consider special cases.
When $\calk$ is the nonnegative orthant, $\calz(\calk)$ consists of square matrices with nonpositive off-diagonal entries, 
$\pi(\calk)$ consists of nonnegative matrices, and LL$(\calk)$ consists of diagonal matrices. Consequently, proper polyhedral cones can be handled via isomorphism arguments. Moving away from proper polyhedral cones,  
in this paper, we focus on the second order cone (also called the Lorentz cone or the ice-cream cone) in the Hilbert space $\rn$, $n>1$, defined by:
\begin{equation}\label{soc}
\call:=\{(t,u)\tp:t\in \mathbb{R},\,u\in \mathbb{R}^{n-1},\,t\geq ||u||\}.
\end{equation}
This cone, being an example of a symmetric cone, appears prominently in conic optimization 
\cite{AlizadehGoldfarb2003}. For this cone,
 Stern and Wolkowicz \cite{SternWolkowicz1991} have shown that {\it $A\in \calz(\call)$ if and only if for some real number $\gamma$, the matrix $\gamma\,J-(JA+A\tp J)$ is positive semidefinite,} where 
$J$ is the diagonal matrix $\diag(1,-1,-1,\ldots, -1)$. 
Another result of Stern and Wolkowicz (\cite{SternWolkowicz1994}, Theorem 4.2) asserts that  
\begin{equation}\label{stern-wolkowicz1994}
\calz(\call)= \mbox{LL}(\call)-\pi(\call).
\end{equation}
(Going in the reverse direction, in  a recent paper, Kuzma et al., \cite{Kuzma2015} have shown that
for an irreducible symmetric cone $\calk$, the equality
$\calz(\calk)= \mbox{LL}(\calk)-\pi(\calk)$ holds only when $\calk$ is isomorphic to $\call.$)
\\
Characterizations of $\pi(\call)$ and LL$(\call)$ appear, respectively, in \cite{LoewySchneider1975} and \cite{TaoGowda2013}.
\\

In this paper, we  describe $\calz(\call)$ and  its interior, boundary, and dual  
in terms of the semidefinite cone and the so-called  copositive and completely positive cones induced by $\call$ 
(or its boundary $\partial(\call)$) see below for the definitions. In particular, we describe the 
dual of $\calz(\call)$ as a 
slice of the semidefinite cone and  also of the  completely positive cone of $\call$. This provides an example of an 
instance where a conic linear optimization problem over a completely positive cone is reduced to a semidefinite problem. 
To elaborate,  consider  $\rn$, 
the Euclidean $n$-space of (column) vectors  with the usual inner product, $\mathbb{R}^{n\times n}$, 
the space of all real $n\times n$  matrices with the inner product $\langle X,Y\rangle =\tr(X^TY)$, 
and $\Sn$, the subspace of all 
real $n\times n$ symmetric matrices in $\mathbb{R}^{n\times n}$. Corresponding to  a closed cone $\calc$ 
(which is not necessarily convex) in $\mathbb{R}^n$, let 
$$\E_\calc:=\copos(\calc):=\Big\{A\in \Sn: x\tp Ax\geq 0,\,\,\forall\,\,x\in \calc\Big\}$$
denote the {\it copositive cone of $\calc$} and 
$$\calk_\calc:=\compos(\calc):=\Big\{\sum uu\tp:u\in \calc\Big \}$$
denote the {\it completely positive cone of $\calc$},  where 
the sum is a finite sum of objects. 
When $\calc=\mathbb{R}^n$, these two cones coincide with the semidefinite cone $\Snp$;
when $\calc=\mathbb{R}^n_+$, these reduce, respectively, to the (standard) copositive cone and completely positive cone. 
All these cones appear prominently in conic optimization. 
A result of Burer \cite{Burer2012} (see also, \cite{Burer2009,Dickinson2013}) 
says that any nonconvex quadratic programming problem over a closed cone with additional linear and binary constraints 
can be reformulated as a linear program over a suitable completely positive 
cone. For this and other reasons, there is a strong interest in understanding copositive and completely positive cones. 
For the closed convex cones $\E_\calc$ and $\calk_\calc$, various structural properties (such as the interior, boundary) 
 as well  as duality, irreducibility,  and homogeneity properties, have been investigated in the literature,  
see for example,  \cite{SturmZhang2003,Dickinson2011,DurStill2008,GowdaSznajder2013}.  
Taking  $\calc$ to be one of $\rn$, $\call$, or $\partial(\call)$, we show that 
\begin{equation} \label{dual of zl}
\calz(\call)^*=\{B\in \mathbb{R}^{n\times n}: \langle B,I\rangle=0, -JB\in \calk_\calc\}
\end{equation}
and deduce the  equality of slices    
\begin{equation}\label{slice relation}
\{X\in\mathbb{R}^{n\times n}:\langle J,X\rangle=0,\, X\in \Snp\}=\{X\in\mathbb{R}^{n\times n}: \langle J,X\rangle=0,\, X\in 
\calk_\calc\}.
\end{equation}

\section{Preliminaries}
In a (finite dimensional real) Hilbert space $(\calh,\lng\cdot,\cdot\rng)$, a nonempty set $\calk$ is said to be a {\it closed convex cone} if it
is closed and $tx+sy\in \calk$ whenever $x,y\in \calk$ and $t,s\geq 0$ in $\re$. Such a cone is said to be {\it proper} if $\calk\cap
-\calk=\{0\}$ and has nonempty interior. Corresponding to a closed convex cone $\calk$, we define its dual in $\calh$ as the set 
\[\calk^*=\{x\in\calh:\lng x,y\rng\ge0,\textrm{ }\forall y\in\calk\}\]
and the {\it complementarity set of $\calk$} as the set 
$\{(x,y): x\in \calk,\,y\in \calk^*,\,\langle x,y\rangle=0\}.$
We say that a linear transformation $A:\calh\rightarrow \calh$ is {\it copositive} on $\calk$ if $\langle Ax,x\rangle \geq 0$ for all $x\in \calk$. We also let 
$\pi(\calk)=\{A: A(\calk)\subseteq \calk\},$ where $A$ denotes a linear transformation on $\calh$.
For a set $S$ in $\calh$, we denote the closure, interior, and the boundary by
$\overline{S}$, $S^\circ$, and $\partial(S)$ respectively.
Throughout this paper, we use the summation sign $\sum $ to describe a finite sum of  objects.

We will be considering closed convex cones in the space $\calh=\rn$  which 
carries the usual inner product and in the space $\re^{n\times n}$ which carries the inner product
$\lng X,Y\rng:=\tr(X\tp Y),$ where the trace of a square matrix is the sum of its diagonal entries.
In $\nn$, $\Sn$ denotes the subspace of all symmetric matrices and $\cala^n$ denotes the subspace of all skew-symmetric matrices. We note that $\nn$ is the orthogonal direct sum of $\Sn$ and  $\cala^n$.


We recall some (easily verifiable) properties of the second order cone $\call$  given by (\ref{soc}).  $\call$ is a self-dual cone in $\rn$, that is, $\call^*=\call$; its interior and boundary are given, respectively, by
$$\call^\circ=\Big \{(t,u)\tp:\,t>||u||\Big \},$$  
$$\partial (\call)=\Big \{(t,u)\tp:\,t=||u||\Big \}=\Big \{\alpha\,(1,u)\tp: \alpha\geq 0,\,||u||=1\Big \}.$$
We also have 
\begin{equation} \label{complementary pair}
\Big [0\neq x,y\in \call,\,\langle x,y\rangle =0\Big ]\Rightarrow x=\alpha\,(1,u)\tp\,\,\mbox{and}\,\,y=\beta\,(1,-u)\tp, \,\,\mbox{for some}\,\,\alpha,\beta>0\,\,\mbox{and}\,\,||u||=1.
\end{equation}

\gap

For a closed cone $\calc$ in $\rn$, we consider the copositive cone $\E_\calc$ and the completely positive cone 
$\calk_\calc$ (defined in the Introduction). Note that these are cones of symmetric matrices. \\

{\it In the Hilbert space $\Sn$ (which carries the inner product from $\nn$), the following  hold.
\begin{itemize}
 \item [$(1)$] $\calk_\calc$ is the dual cone of $\E_\calc$ \cite{SturmZhang2003}. 
\item [$(2)$] When $\calc-\calc=\rn$, both $\E_\calc$ and $\calk_\calc$ are  proper cones (\cite{Gowda2012}, Proposition 2.2). In particular, this holds when $\calc$ is one of $\rn$, $\call$, or $\partial(\call)$. 
\item [$(3)$] We have
$\E_{\rn}=\Snp\subset \E_{\call}\subset \E_{\partial(\call)},$ or equivalently, 
$\calk_{\partial(\call)}\subset \calk_{\call}\subset \calk_{\rn}=\Snp.$

\end{itemize}
}

\section{Main results}

In this section,  we provide a  closure-free description of $\calz(\call)$ and, additionally, describe the dual, interior, and the boundary of $\calz(\call)$.
We recall that $J=\diag(1,-1,-1,\ldots, -1)$ and $\an$ denotes the set of all skew-symmetric matrices in $\nn$.  

\begin{theorem}\label{tzl}
{\it 
Let $\calc$ denote one of $\rn$, $\call$, or $\partial(\call)$. Then,
\begin{equation}\label{zl without closure}
	\calz(\call)=\re\,\I-J(\E_\calc+\cala^n).
\end{equation}
}
\end{theorem}

\gap

\begin{proof}
Let $A\in \calz(\call).$ From the result of Stern and Wolkowicz \cite{SternWolkowicz1991} mentioned in the Introduction, 
we have 
$$2\gamma J-(JA+A\tp J)=2P$$ for some
	$\gamma\in\re$ and $P\in\snp$. Hence, $JA+(JA)\tp=2(\gamma J-P)$, which implies 
	\begin{equation}\label{eprej}
		2JA=JA+(JA)\tp-\lf[(JA)\tp-JA\rg]=2(\gamma J-P)-2Q,
	\end{equation}
	where $2Q=(JA)\tp -JA$ is skew-symmetric.
	Since $J^2=I$, this leads to  
$$A=\gamma\,\I -J(P+Q),$$
 where $P\in \Snp$ and $Q\in \cala^n$.
As $\Snp\subset \E_\call\subset \E_{\partial(\call)}$, this proves  that 
\begin{equation}\label{inclusions}
	\calz(\call)\,\subseteq\, \re\,\I-J(\Snp+\cala^n)\,\subseteq\, \re\,\I-J(\E_\call+\cala^n)\,\subseteq\,
	\re\,\I-J(\E_{\partial(\call)}+\cala^n).
\end{equation}
Now, to see the reverse inclusions,  
suppose $A=\gamma\,\I-J(P+Q)$ for some
$\gamma\in \re$, $P\in \E_{\partial(\call)}$, and $Q$  skew-symmetric. Let $0\neq x,y\in \call$ with $\langle x,y\rangle=0$. 
By (\ref{complementary pair}), $x$ and $y$ are in $\partial(\call)$, and $Jy$  is a positive multiple of $x$. 
Hence, $\langle Px,Jy\rangle \geq 0$ as $P\in \E_{\partial(\call)}$ and $\langle Qx,Jy\rangle=0$ as $Q$ is skew-symmetric. Thus, 
$$\langle Ax,y\rangle=\gamma\langle x,y\rangle-\langle JPx,y\rangle+
\langle JQx,y\rangle=-\langle Px,Jy\rangle+\langle Qx,Jy\rangle\leq 0.$$ 
This shows that $A\in \calz(\call)$ and so, inclusions in (\ref{inclusions}) turn into equalities. Thus we have (\ref{zl without closure}).
\end{proof}

\noindent{\bf Remarks.}
From the above theorem, we have
        $$\re\,\I-J(\Snp+\cala^n)\,=\, \re\,\I-J(\E_\call+\cala^n)\,=\,
        \re\,\I-J(\E_{\partial(\call)}+\cala^n).
$$
Multiplying throughout by $J$ and noting $-\an=\an$, we get the equality of sets
$$(\re\,J-\Snp)+\cala^n\,=\,(\re\,J-\E_\call)+\cala^n\,=\,
        (\re\,J-\E_{\partial(\call)})+\cala^n,
$$
where each set is a sum of $\an$ and a subset of $\sn$. Since $\nn=\sn+\an$ is an (orthogonal) direct sum decomposition, 
we see that
\begin{equation}\label{equality of sets}
        \re\,J-\Snp\,= \,\re\,J-\E_\call\, =\, \re\,J-\E_{\partial(\call)}.
\end{equation}
These equalities can also be established via different arguments.
A result of 
Loewy and Schneider \cite{LoewySchneider1975} asserts that 
{\it   A symmetric matrix $X$ is  copositive on $\call$ if and only if
 there exists $\mu\geq 0$ such that $X-\mu\,J\in \Snp$.} (This is essentially a consequence of the so-called S-Lemma \cite{PolikTerlaky2007}: 
If $A$ and $B$ are two symmetric matrices with $\langle Ax_0,x_0\rangle >0$ for some $x_0$ and $\langle Ax,x\rangle \geq 0\Rightarrow \langle Bx,x\rangle \geq 0$, 
then there exists $\mu\geq 0$ such that $B-\mu A$ is positive semidefinite.) 
This result gives the equality
$$\E_\call=\snp+\re_+\,J$$
and consequently $\re\,J-\snp=\re\,J-\E_\call$.   
The equality
$$
\E_{\partial(\call)}=\snp+\re\,J
$$
can be seen via an application of Finsler' theorem \cite{PolikTerlaky2007} that says that if $A$ and $B$ are two symmetric matrices with $\left [x\neq 0, \langle Ax,x\rangle=0\right ] \Rightarrow \langle Bx,x\rangle > 0,$ 
then there exists $\mu\in \re$ such that $B+\mu A$ is positive semidefinite. (For $M\in \E_{\partial(\call)}$ and vectors $u,v\in \call^\circ$, one has $\langle Jx,x\rangle =0\Rightarrow \langle M_k x,x\rangle >0$, where $k$ is a natural number and 
$M_k:=M+\frac{1}{k}uv\tp$. When $M_k+\mu_k\,J$ is positive semidefinite, it follows that the sequence $\mu_k$ is bounded.) From this equality, one gets  
$\re\,J-\snp=\re\,J-\E_{\partial(\call)}$. 
\\

Our next result deals with the dual of $\calz(\call)$.

\begin{theorem} \label{dual of zl-2} {\it Let $\calc$ denote one of $\rn$, $\call$, or $\partial(\call)$. Then,   
$$\calz(\call)^*=\{B\in \mathbb{R}^{n\times n}: \langle B,I\rangle=0, -JB\in \calk_\calc\}.$$
In particular, (\ref{slice relation}) holds.
}
\end{theorem}

\gap

\begin{proof}
We fix $\calc$.  From (\ref{zl without closure}), we see that 
$B\in \calz(\call)^*$ if and only if  
$$0\leq \langle B, \gamma\,\I-J(P+Q)\rangle$$
for all $\gamma$ real, $P$ in $\E_\calc$, and $Q$  in $\cala^n$. Clearly, this holds if and only if 
$$\langle B,I\rangle =0,\,\langle -JB,P\rangle \geq 0,\,\mbox{and}\,\,\langle -JB,Q\rangle =0$$
for all $\gamma$, $P$, and $Q$ specified above. Now, with the observation that a (real) matrix is orthogonal to all skew-symmetric matrices  in $\nn$ if and only if it is symmetric,  this  further simplifies to
$$\langle B,I\rangle =0\,\,\mbox{and}\,\,-JB\in \E_\calc^*,$$
where $\E_\calc^*$ is the dual of $\E_\calc$ computed in $\Sn$. Since $\calk_\calc=\E_\calc^*$ in $\Sn$, we see that 
$B\in \calz(\call)^*$ if and only if $\langle B,I\rangle =0\,\,\mbox{and}\,\,-JB\in \calk_\calc$. This completes the proof.
\end{proof}

We remark that (\ref{slice relation}) can be deduced directly from (\ref{equality of sets}) by taking the duals in $\sn$.
\\

In our final result, we describe the interior and boundary of  $\calz(\call)^\circ$.
First, we recall some definitions from \cite{FanTaoRavindran2017}. 
Let 
 $$\Omega:=\Big \{(x,y)\in \call\times \call: ||x||=1=||y||\,\,\mbox{and}\,\, \langle x,y\rangle=0\Big \}.$$ 
It is easy to see that $\Omega$ is compact and, from (\ref{complementary pair}), 
\begin{equation}\label{another form of Omega}
\Omega=\Big \{(x,Jx):x\in \partial(\call),\,||x||=1\Big \}.
\end{equation}
For any $A\in \re^{n\times n}$, let 
$$\gamma(A):=\max\Big \{\langle Ax,y\rangle:(x,y)\in \Omega\Big \}.$$
 Note that $A\in \calz(\call)$ if and only if $\gamma(A)\leq 0$. 
We say that $A\in \re^{n\times n}$ is a {\it 
strict-$\calz$-transformation on $\call$} if 
$$\Big [0\neq x,y\in \call,\,\langle x,y\rangle=0\Big ]\Rightarrow \langle Ax,y\rangle<0.$$
The set of all such transformations is denoted by $str(\calz(\call))$. 
For $A\in \nn$, the following statements are shown in \cite{FanTaoRavindran2017}, Theorem 3.1:
$$\gamma(A)<0\Longleftrightarrow A\in \calz(\call)^\circ \Longleftrightarrow A\in str(\calz(\call))$$
and
$$\gamma(A)=0\Longleftrightarrow A\in \partial(\calz(\call)).$$ 

Recall that $\E_\call$ consists of all symmetric matrices that are copositive on $\call$. We say that a symmetric matrix 
$P$ is {\it strictly copositive on $\call$} if
$0\neq x\in \call\Rightarrow \langle Px,x\rangle>0$;  the set of all such matrices is denoted by $str(\E_\call)$.
Similarly, one defines $str(\E_{\partial(\call)})$.  

\begin{corollary} 
{\it The following statements hold:  
$$\calz(\call)^\circ=str(\calz(\call))=\re\,\I-J\,\Big (str(\E_{\partial(\call)})+\an\Big )$$
and 
$$\partial(\calz(\call))=\re\,\I-J\,\Big (\partial_{*}(\E_{\partial(\call)})+\an\Big ),$$
where $\partial_{*}(\E_{\partial(\call)})$ denotes the boundary of $\E_{\partial(\call)}$ in $\sn$.
}
\end{corollary}

\gap

\begin{proof}
We first deal with the interior of $\calz(\call)$.  The equality 
$$\Big \{A\in \re^{n\times n}:\gamma(A)<0\Big \}=\calz(\call)^\circ=str(\calz(\call))$$ has already been observed in  
\cite{FanTaoRavindran2017}, Theorem 3.1.
To see the first assertion, we show that $\gamma(A)<0$ if and  only if $A=\theta\, I-J(P+Q)$ for some $\theta\in \re$, 
$P$ (symmetric) strictly copositive on $\partial(\call)$, and $Q$ skew-symmetric. 
Suppose $\gamma(A)<0$. Then, for any $\theta\in \re$, 
$$\max\Big \{\big \langle (A-\theta\,I)x,y\big \rangle: (x,y)\in \Omega\Big \}<0,$$
which, from  (\ref{another form of Omega}) becomes 
$$\min\Big \{\big \langle J(\theta\,I-A)x,x\big \rangle:\,x\in \partial(\call), ||x||=1\Big \}>0.$$
Now, fix $\theta$ and let $J(\theta\,\I-A)=P+Q$, where $P\in \sn$ and $Q\in \an$. As $\langle Qx,x\rangle=0$ for any $x$, the above inequality implies that $\min\Big \{\big \langle Px,x\big \rangle:\,x\in \partial(\call), ||x||=1\Big \}>0.$ This proves that 
 $P$ is strictly
copositive on $\partial(\call)$. Rewriting $J(\theta\,\I-A)=P+Q$, we see that $A=\theta\,\I-J(P+Q)$ which is of the required form.
\\
To see the converse, suppose $A=\theta\,\I-J(P+Q)$, where $\theta\in \re$,
$P$ (symmetric) strictly copositive on $\partial(\call)$, and $Q$ skew-symmetric.
Using (\ref{another form of Omega}), we can  easily verify that $\gamma(A)<0$. 
 Thus, $A\in str(\calz(\call))$.   
\\
An argument similar to the above will show that 
$\gamma(A)=0$ if and  only if 
$A=\theta\, I-J(P+Q)$ for some $\theta\in \re$,
$P\in \partial_{*}(\E_{\partial(\call)})$, and $Q$ skew-symmetric. 
This gives the statement regarding the  boundary of $\calz(\call)$.
\end{proof}

We end the paper with a  remark dealing with conic linear programs. Motivated by the result of Burer (mentioned in the Introduction), we consider a conic linear program on a completely positive cone $\calk_\calc$ (where $\calc$ is a closed cone):
$$\min\Big \{\langle c,x\rangle: Ax=b, x\in \calk_\calc\Big \}.$$
While such a problem is generally hard to solve, we ask:
(When) can we replace $\calk_\calc$ by $\Snp$ and thus reduce the above problem  to the semidefinite programming problem
$\min\Big \{\langle c,x\rangle: Ax=b, x\in \Snp\Big \}?$
Just replacing $\calk_\calc$ by $\snp$ without handling the constraint $Ax=b$ is not viable as $\calk_\calc=\Snp$ if and only if $\calc\cup -\calc=\rn$ (which fails to hold when $n>1$ and $\calc$  is pointed), see \cite{GowdaSznajder2013}. 
While we do not answer this broad question, we point out, as a consequence of (\ref{slice relation}) that for any $C\in \Sn$, 
$$\min\Big \{\langle C,X\rangle: \langle X,J\rangle=0, X\in \calk_\call\Big \}=\min\Big \{\langle C,X\rangle: \langle X,J\rangle=0, 
X\in \snp\Big \}.$$ 


\bibliographystyle{unsrt}
\bibliography{zstar}
                                                                            
\end{document}